\documentclass{amsart}
\usepackage{amssymb, graphics, amsmath, epsfig, changebar}

\theoremstyle{plain}
\newtheorem{theorem}{Theorem}
\newtheorem{corolary}{Corolary}

\newtheorem{lemma}{Lemma}

\newtheorem{remark}{Remark}

\theoremstyle{definition}
\newtheorem{definition}{Definition}

\theoremstyle{remark}

\numberwithin{equation}{section}
                  
\begin{document}

\title[Diagrammatic unknotting of knots and links in $\mathbb RP^3$]{Diagrammatic unknotting of knots and links in the projective space}
\author{Maciej Mroczkowski}

\address{Department of Mathematics,
         Uppsala University, 751 06 Uppsala, Sweden}
\email{mroczkow@math.uu.se}        
\begin{abstract}
In the classical knot theory there is a well-known notion of descending diagram. From an arbitrary diagram one can easily obtain, by some crossing changes, a descending diagram which is a diagram of the unknot or unlink.
In this paper the notion of descending diagram for knots and links in $\mathbb R^3$ is extended to the case of nonoriented knots and links in the projective space. It is also shown that this notion cannot be extended to oriented links.
\end{abstract}
\maketitle

\section{introduction}
Embeddings are dense in the space of immersions of a curve to a
3-manifold. Hence any immersion of a collection of circles to a
3-manifold can be turned by a small regular homotopy into a
differentiable embedding. Usually, the ambient isotopy type of an
embedding, which can be obtained from a given immersion by an
arbitrarily small (in $C^1$-topology) regular homotopy, is not entirely
determined by the immersion. On the other hand, for any immersion there
are ambient isotopy types which cannot be obtained from it. 

The main results of this paper imply the following theorem.

\begin{theorem}
Let $S$ be a smooth closed 1-manifold. For any immersion
$f:S\to\mathbb
RP^2$, its composition with the standard inclusion $in: \mathbb RP^2\to\mathbb RP^3$ is regularly homotopic via an arbitrarily small regular homotopy to an embedding $g:S\to\mathbb RP^3$, which depends,
up to ambient isotopy and composition with a  self-diffeomorphism of $S$,
only on the homotopy class of $f$.
In other words, $g(S)$ is ambiently isotopic to a standard nonoriented unlink $L_{p,q}\subset\mathbb RP^3$ which depends, up to ambient isotopy, only on the number $p$ of its components contractible in $\mathbb RP^3$ and the number $q$ of its components non contractible in $\mathbb RP^3$.
\end{theorem}

From this theorem one cannot eliminate  self-diffeomorphisms of $S$.
This is also proven below (see Section \ref{counter}).

This paper presents a way to unknot knots and unlink links in the real projective 3-space $\mathbb RP^3$ and the results formulated above appear as straightforward corollaries. It is shown how to obtain a diagram of the unknot or unlink, starting from an arbitrary diagram, and performing some crossing changes on it. This is done through an extension of the classical notion of descending diagram to diagrams of knots and links in $\mathbb RP^3$.

In the case of knots in $\mathbb R^3$, the notion of descending diagram was used to study some knot invariants such as Jones polynomial, Homfly polynomial or finite type invariants.
One can calculate these invariants using appropriate skein relations and the fact that it is possible to make any link diagram descending, which would turn it to a diagram of the unlink. Descending diagrams were used to define the Homfly polynomial (see for instance \cite{LM}).

The notion of descending diagram for knots in $\mathbb RP^3$ can be used to study some invariants of these knots. It can be useful when considering the Homflypt module of $\mathbb RP^3$.
The Jones polynomial was extended in \cite{JD} by J. V. Drobotukhina to knots in $\mathbb RP^3$. An algorithm that makes a diagram descending gives an alternative way to calculate this polynomial.

I wish to thank Oleg Viro for stimulating conversations and for his help.

\section{Unknotting knot and link diagrams in $\mathbb R^3$}\label{r3descending}

In the case of knots in $\mathbb R^3$ there is a well-known way to obtain the unknot by doing some crossing changes on a diagram.
We choose a starting point and proceed from this point according to some orientation of the knot. When we meet a crossing for the first time and arrive at the lower branch of the crossing, we change it in order to make this branch upper.
If we arrive at a crossing for a second time, we leave it unchanged. Finally we get back at the initial point.
When we have done all these changes the diagram takes a special form: it is descending.

We can imagine that the result is a diagram of a knot in which we descend from the initial point all the way (a $z$ coordinate is decreased if $z$ is the axis along which we project the knot to the diagram). When we get back at a point which has the same projection as the initial point (but is below it) we join these two points with a segment. The resulting knot is the unknot.

More generally, a link diagram can be unlinked with appropriate crossing changes, by putting different components one above an other (choosing some order) and making each component descending as above.

\section {knot and link diagrams in $\mathbb RP^3$}
\subsection{Link diagrams}
\begin{figure}[h]
\scalebox{1}{\includegraphics{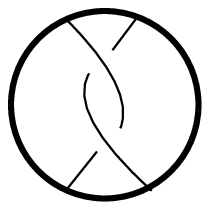}}
\caption{}
\label{knot1}
\end{figure}

A link diagram in $\mathbb RP^3$ is a disk with a collection of generically immersed arcs.
An arc is a compact connected 1-manifold with or without boundary.
The endpoints of arcs with boundary are on the boundary of the disk, divided into pairs of antipodal points. Each double point of the immersions or {\it crossing} of the diagram is endowed with information of over- and undercrossing. An example of a knot diagram in $\mathbb RP^3$ is shown in Figure \ref{knot1}.

A link diagram is constructed from a link $L$ in the following way:

$\mathbb RP^3$ is represented as a ball $D^3$ with antipodal points of the bounding sphere identified. The link $L$ in $\mathbb RP^3$ is lifted to $L'$ in $D^3$. We can suppose that the poles of the ball are not in $L'$. Let $p$ be the projection of $L'$ to the equatorial disk $D^2$ where a point in $L'$ is projected along the metric circle in $D^3$ passing through this point and the poles of the ball.

We assume that $L$ satisfies the following conditions of general position: $p(L')$ contains no cusps, points of tangency, triple points; $L'$ intersects transversally the boundary of the ball; no two points in $L'$ lie on the same arc of the great circle joining the poles of the ball in the boundary of the ball.

The information of over- and undercrossings comes from some orientation of the circles along which $L'$ is projected to $D^2$ (for instance from north to south).

If the link $L$ is oriented, we get naturally an {\it oriented} link diagram for which each arc is oriented. An orientation of a link diagram gives rise to a cyclic ordering of arcs (when we travel on $L$ according to the orientation, we meet the arcs in this order).

\begin{figure}[h]
\scalebox{1}{\includegraphics{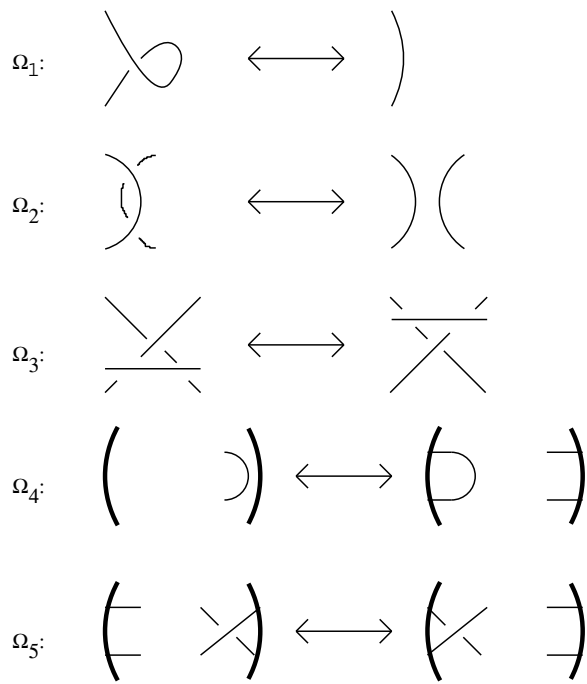}}
\caption{Reidemeister moves}
\label{reidemeister}
\end{figure}

As for diagrams of links in $\mathbb R^3$, there are Reidemeister moves for diagrams of links in $\mathbb RP^3$. The five of them are pictured in Figure \ref{reidemeister}. These moves appeared in \cite{JD}.

\subsection{Nets, diagrammatic components, arc distance and first pass}
A {\it net} is the projective plane $\mathbb RP^2$ with a distinguished line, called {\it the line at infinity}, and a collection of generically immersed circles together with information of over- and undercrossing for each double point.
We can associate to each diagram $D$ of a link {\it its net} obtained from $D$ by identifying the antipodal points of the boundary circle of $D$, with the line at infinity coming from this boundary circle.

Let $D$ be a diagram of a link and $b$ the part of $D$ coming from a component, $L_b$, of the link. We will say that $b$ is a {\it diagrammatic component} of $D$. Suppose that $L_b$ is oriented, then $b$ is also oriented.
Let $P$ and $Q$ be two points in the interior of some arcs of $b$. Then the {\it arc distance} between $P$ and $Q$ is defined to be the number of times the line at infinity is crossed in the net of $D$, if one travels from the image of $P$ to the image of $Q$ in the net, according to the orientation of the image of $b$ in the net.
Suppose that $X$ is a crossing of $D$ such that at least one of its branches is in $b$. Then the {\it first pass} of $X$ from $P$ is, by definition, the branch of $X$ whose image in the net of $D$ is passed first, if one travels from the image of $P$ in the net, according to the orientation of the image of $b$ in the net.

\subsection{Unknots in $\mathbb RP^3$}
The fundamental group of $\mathbb RP^3$ has two elements. In each of them there is a simple loop that is naturally called {\it unknot}. A planar circle and a projective line are two unknots, up to isotopy. They are the only knots for which there are diagrams without crossings.  A knot in $\mathbb RP^3$ is homotopic to one of the unknots depending on the element of the fundamental group it realizes. Thus we can deform it to an unknot by a sequence of isotopies, which correspond to planar isotopies and Reidemeister moves on the level of diagrams, and some homotopies, namely the ones which correspond to crossing changes on the level of diagrams.

In the next section we will see that for any diagram of a knot we can obtain a diagram of an unknot solely by some crossing changes.

\section{Unknotting knot diagrams in $\mathbb RP^3$}

A natural question is whether, for knots in $\mathbb RP^3$, there is a way to obtain a diagram of the unknot, by changing some crossings on an arbitrary diagram. The goal of this section is to give a positive answer to this question.

A {\it basepoint} is a distinguished point of a diagram, distinct from crossings and endpoints of arcs. A diagram with a basepoint is called {\it based diagram}.

\begin{definition}
A based oriented diagram $D$ is called {\it descending} provided that for every crossing $X$ of $D$, the first pass of $X$ from the basepoint is an overpass (resp. underpass), if the arc distance between the basepoint and this first pass is even (resp. odd).
\end{definition}

\begin{theorem}\label{theorem_knots}
Let $D$ be a based oriented diagram of a knot in $\mathbb RP^3$. If $D$ is descending, then $D$ is a diagram of an unknot.
\begin{proof}
Suppose that $D$ is descending.

First, note that if $D$ consists of a single arc without boundary, then $D$ is a diagram of 0-homologous unknot, because it is descending in the classical sense.

Now suppose that $D$ is not of that type. Denote the arcs of $D$ by $a_1, ..., a_n$, $n\ge 1$, where $a_1$ contains the basepoint and $a_2, ..., a_n$ are ordered according to the orientation of $D$.

The arc $a_1$ can be divided in two parts: one that comes after the basepoint (according to the orientation) and the other one that comes before the basepoint. Denote the first one by $a_1^a$ and the second one by $a_1^b$.

An arc, or a part of it is said to be {\it below} another one if at each crossing between the two of them the branch of the first one is below the branch of the second one. If $a$ is below $b$, we write $a\le b$.
It is easy to see that in the descending diagram $D$ the following relation holds:

$$a_2\le a_4\le a_6\; ...\le a_1^b\; ...\le a_5\le a_3\le a_1^a$$

Observe that each arc is descending or ascending. Also $a_1^a$ and $a_1^b$ are descending or ascending.

We will show that $D$ is a diagram of an unknot by constructing a sequence of Reidemeister moves from $D$ to a diagram without crossings.  
As each arc is descending or ascending we kill all the crossings between an arc and itself with some $\Omega_1-\Omega_3$ moves. We do the same with $a_1^a$ and $a_1^b$.

Now consider Figure \ref{isotopy}. We want to reduce the number of arcs by eliminating $a_2$. One of its boundary points, say $P$, is antipodal to a boundary point of $a_1^a$. Denote the other boundary point of $a_2$ by $Q$. Denote by $P'$, resp. $Q'$ the antipodal points of $P$ resp. $Q$.

\begin{figure}[h]
\scalebox{1}{\includegraphics{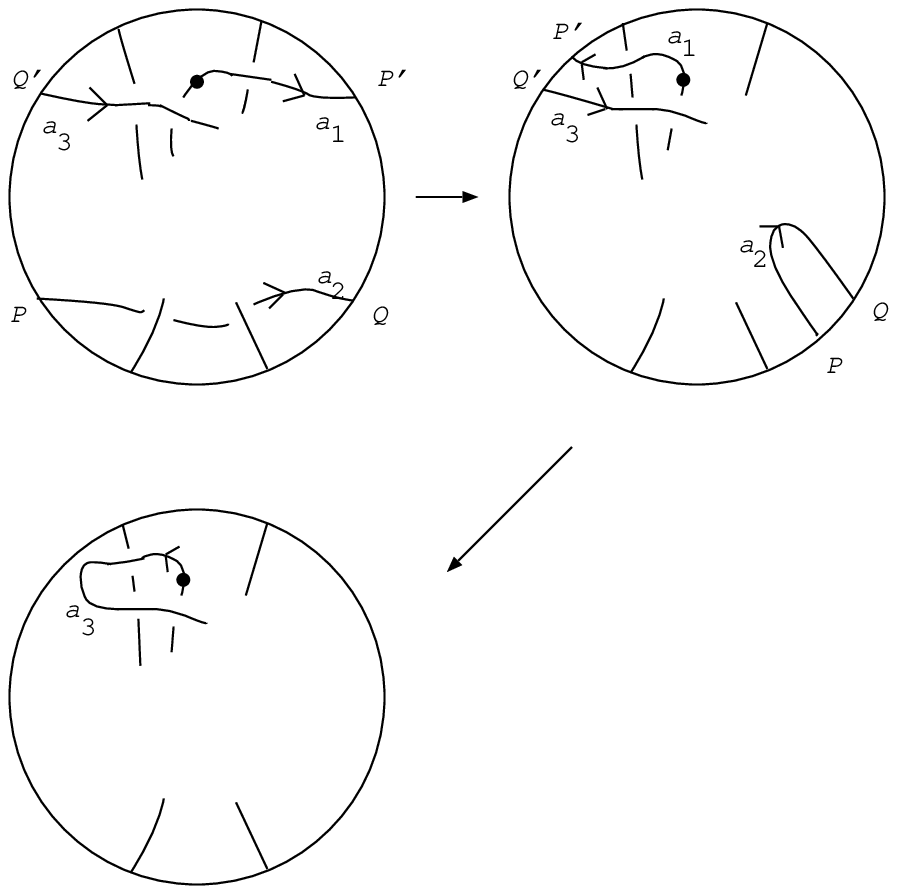}}
\caption{}
\label{isotopy}
\end{figure}

In order to eliminate $a_2$ by $\Omega_4$ move, we want to have no endpoints between the two endpoints of $a_2$, on the boundary of the disk. For this purpose we move $P$ towards $Q$, below some endpoints. At the same time $P'$ moves towards $Q'$ above some other endpoints. This corresponds to several applications of $\Omega_5$ move.
Secondly, we perform several $\Omega_1-\Omega_3$ moves to kill all the crossings between $a_2$ and other arcs. It is possible because $a_2$ is below all other arcs.
Finally, we can perform $\Omega_4$ move on $a_2$. Now instead of $a_1^a$, $a_2$ and $a_3$ we have a single part of arc, call it $a_3^a$. It is unknotted and above all other arcs. There may be some self-crossings of $a_3^a$ and in that case we kill them with $\Omega_1-\Omega_3$ moves. $a_1^b$ is renamed $a_3^b$.
The arcs are now positioned in the following way:
$$a_4\le a_6\; ...\le a_3^b\; ...\le a_5\le a_3^a$$

We can repeat the process with $a_3^a$, $a_4$ and $a_5$. As long as there are at least three arcs or parts of arcs left this is possible. Finally we will end up with a single arc. It will be unknotted and we can kill all crossings again. Thus it will be 0-homologous or non zero-homologous unknot.
\end{proof}
\end{theorem}

\begin{corolary}
Let $D$ be a diagram of a knot in $\mathbb RP^3$. By making some crossing changes on it, we can obtain a diagram of the unknot.
\begin{proof}
Choose orientation and basepoint on $D$. Let $D'$ be the diagram obtained from $D$ by making the necessary crossing changes that make $D'$ descending. Then $D'$ is a diagram of an unknot by Theorem \ref{theorem_knots}.
\end{proof}
\end{corolary}

\begin{remark}
In the proof of Theorem \ref{theorem_knots}, there are two types of simplifications of the diagram. The first one is to make the arcs unknotted. In this case we use only $\Omega_1$, $\Omega_2$ and $\Omega_3$ moves. The second is to eliminate some arcs. In this case we use only $\Omega_2$ and $\Omega_3$ moves that involve at least two different arcs as well as $\Omega_4$ and $\Omega_5$ moves.
In the proof we alternated these two types of simplifications.

Another proof is possible in which any descending diagram is changed to a diagram with no crossings by a sequence of simplifications of the second type (the number of arcs will be reduced to one), followed at the end by a simplification of the first type.
\end{remark}

\begin{remark}
In the case of a diagram that is descending and which represents 0-homologous unknot, consider the diagram with a single arc obtained after a sequence of simplifications of the second type defined in the previous remark. This diagram will be descending in the classical sense (when we move an ascending arc through the boundary it becomes descending and vice versa). In this way we see that the notion of descending diagram in $\mathbb RP^3$ is similar to the same notion in $\mathbb R^3$ even for some diagrams of non affine knots.
\end {remark}

\section{unlinks in $\mathbb RP^3$}

There is no natural notion of unlink in $\mathbb RP^3$. A link has two types of components: the 0-homologous and the non 0-homologous ones (they will be called {\it 1-homologous}).
For the first type there is no problem to see how they should look in an unlink: like in $\mathbb R^3$ they should be unknots and each of them should be in a ball which does not intersect the other components.
For the second type, different definitions of unlink are possible. Note that two 1-homologous components will always intersect in a diagram. Even if we require that an unlink should have a diagram in which any couple of 1-homologous components has a single common crossing and there are no other crossings, there are still many choices for an unlink.
This is related to the configurations of skew lines in \cite{OYV} by O. Ya. Viro and J. V. Drobotukhina.

\begin{figure}[h]
\scalebox{1}{\includegraphics{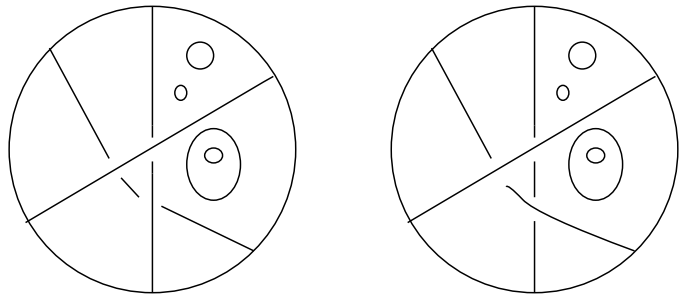}}
\caption{unlinks in $\mathbb RP^3$}
\label{unlinks}
\end{figure}

In Figure \ref{unlinks}, the two links could be taken as unlink with three 1-homologous components and four 0-homologous. But in fact they are not isotopic.
This suggests that one may define a {\it standard unlink}. We consider the following definition of standard unlink:
Take a projective line. Add the next below the first slightly rotating it in the counterclockwise direction. The third will be under the two firsts, also rotated in the same direction. We continue in this way with all projective lines.

An equivalent way to obtain standard unlink is the following:
Take several complex lines in $\mathbb C^2$. They represent projective lines in $\mathbb RP^3$ which form the standard unlink.

For instance in Figure \ref{unlinks} the first is standard unlink, but not the second.

\section{Unlinking link diagrams}
Now the question arises: is it possible to unlink a link diagram (obtain a diagram of a standard unlink) by crossing changes ?

The approach that was used in the case of knots does not extend directly.
The problems arise from the 1-homologous components.

Recall that a diagrammatic component of $D$ is a part of $D$ coming from a component of the link. We will say that a diagrammatic component is {\it 0-homologous} (resp. {\it 1-homologous}), if it comes from a 0-homologous (resp. 1-homologous) component of the link.

\subsection{0-homologous components}
Suppose that we have a diagram of a link, $D$, with $k$ 0-homologous diagrammatic components, $a^1, ..., a^k$, and some 1-homologous diagrammatic components. Choose some basepoint and orientation for each 0-homologous diagrammatic component.

If we make $a^1$ descending, in the same way as it was done for knots (its arcs will alternate below and above the other components), it will become diagram of an unknot, unlinked to the other components. We can continue with $a^2, ..., a^n$ in the same way.

The crossings between the 1-homologous diagrammatic components of $D$ were not changed, so the problem of unlinking a link diagram reduces to the problem of unlinking a diagram with only 1-homologous diagrammatic components.

\subsection{1-homologous components}
\begin{figure}[h]
\scalebox{1}{\includegraphics{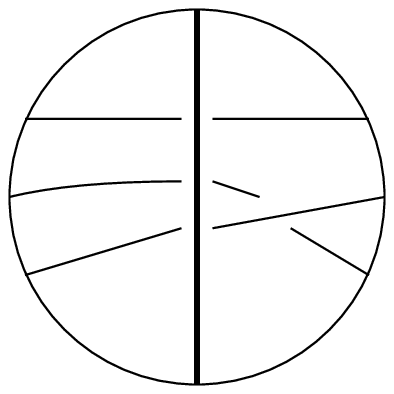}}
\caption{}
\label{bad_1}
\end{figure}

If we try to make the arcs of 1-homologous diagrammatic components of $D$ alternating in the similar way as for 0-homologous diagrammatic components, there will be a problem: after unknotting the first 1-homologous diagrammatic component it will lie below or above everything else and we will not be able to continue with the next one. For instance in Figure \ref{bad_1} we can try to put the component with a single arc below or above the other component. But in either case the resulting diagram will not be a diagram of standard unlink.

Nonetheless, as was stated in the introduction, we will define a notion of descending diagram for arbitrary link diagrams. In order to do this, we need to enhance the data which is used in the definition of descending diagram.

In the next section we will construct the data needed for the 1-homologous diagrammatic components (different from the data needed for the 0-homologous diagrammatic components). For each such component a set of self-crossings will be specified. This set will determine some part of the diagram which will be called dashed part.
We will also specify two antipodal endpoints of arcs on some part of the diagram that is not dashed.

With this data we define the notion of descending diagram. We prove that a descending diagram is a diagram of the standard unlink in two steps. In the first step the dashed part of the diagram is eliminated by a sequence of Reidemeister moves. The remaining diagram has a simpler form and, in the second step, it is transformed into a canonical diagram of standard unlink by another sequence of Reidemeister moves.

\subsection{Data for descending diagrams: simplifying sets.}

Let $D$ be a diagram. Consider a 1-homologous diagrammatic component of $D$, say $b$.
By the set of {\it self-crossings} of $b$, we will mean the set of those crossings of $D$ for which both branches are in $b$.

Let $X$ be a self-crossing of $b$. We can associate to it an orientation of $b$ defined in such a way that, with this orientation, the arc distance between the upper branch of $X$ and the lower one is even. This gives a unique choice of orientation because $b$ is 1-homologous: if we reverse the orientation, the arc distance will be odd.
We call the associated orientation {\it the orientation determined by $X$}.

\begin{figure}[h]
\scalebox{1}{\includegraphics{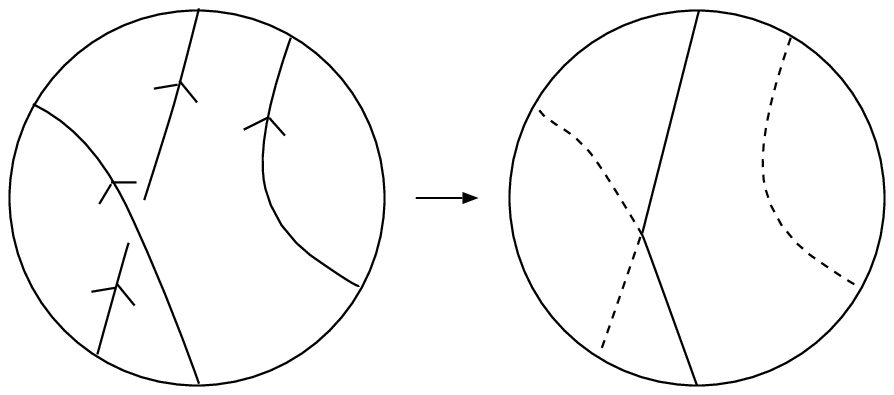}}
\caption{}
\label{dash}
\end{figure}

We will also dash the part of the diagram along which we travel from the upper branch to the lower branch of $X$ according to the orientation determined by $X$. To be precise, we notice that we travel on the net of the diagram, but we can lift it to the diagram itself and dash in this way some part of it. We call it the {\it dashed part determined by $X$}.
An example is shown in Figure \ref{dash}.

A subset $M$ of the set of self-crossings of a 1-homologous diagrammatic component $b$ is called {\it simplifying set} if:

1) Any self-crossing of $b$ is either in $M$ or is such that at least one of its branches is in a dashed part determined by a crossing in $M$.

2) For any two crossings $X, Y$ in $M$, the intersection between the dashed part determined by $X$ and the one determined by $Y$ is empty or consists of some crossings (the dashed parts are disjoint except possibly for some crossings).

To see that it is always possible to construct a simplifying set $M$, we first define a partial order on the set of self-crossings of $b$ in a diagram $D$:

$Y\le X$ if both branches of $Y$ are in the dashed part determined by $X$.
Also $X\le X$ for any $X$. It can be checked easily that this indeed gives a partial order.

A simplifying set $M$ is constructed in the following way: Choose any self-crossing of $b$, say $X_1$, that is maximal with respect to the partial order defined above. It is possible if the set of self-crossings of $b$ is non empty. Put $X_1$ in $M$.
Now suppose that $X_1, ..., X_k$ are already in $M$. Consider the diagram $D'$ obtained from $D$ by removing the dashed parts determined by $X_1, ..., X_k$ (and joining the remaining part of the upper branch to the remaining part of the lower branch for each of the crossings $X_1, ..., X_k$).
Choose any self-crossing of $b$ in $D'$, say $X_{k+1}$, that is maximal. It can be viewed as a crossing in $D$ because the crossings of $D'$ form naturally a subset of the set of crossings of $D$. Put $X_{k+1}$ in $M$. At some point the set of self-crossings of $b$ in $D'$ is empty and $M$ is constructed.

Such set $M$ is indeed a simplifying set:

1) Each crossing is either in $M$ or in a dashed part determined by a crossing in $M$ (otherwise more crossings could be put in $M$). 

2) The intersection between the dashed part determined by $X_1$ and the dashed part determined by any other crossing $Y$ in $M$ is empty or consists of some crossings. This is the case because $X_1$ is maximal and none of the branches of $Y$ is in the dashed part determined by $X_1$.
Similarly, the intersection between the dashed part determined by $X_k$ and the one determined by any crossing $X_l$, with $l>k$, is empty or consists of some crossings. This can be seen by considering the diagram $D'$ obtained from $D$ by removing the dashed parts determined by $X_1, ..., X_{k-1}$.
Thus for any crossings $X, Y$ in $M$, the intersection between the dashed part determined by $X$ and the one determined by $Y$ is empty or consists of some crossings.

\subsection{Simple diagrams}
Diagrams which have only 1-homologous diagrammatic components and which have crossings only between different diagrammatic components will be called {\it simple diagrams}.

Let $D$ be a diagram. Suppose that for each 1-homologous diagrammatic component of $D$, a simplifying set is chosen. Dash the parts of $D$ determined by all crossings in the simplifying sets. Then the diagram obtained from $D$ by removing the dashed parts and all the 0-homologous diagrammatic components, is called the {\it simple diagram} of the data consisting of $D$ and the simplifying sets. This diagram is indeed simple in the sense defined in the preceding paragraph. At the crossings which are in the simplifying sets, we join the remaining part of the upper branch to the remaining part of the lower branch.

\subsection{Simple descending diagrams}\label{simple_desc_dgms}
We will first define the notion of descending diagram for a simple diagram $D$ with a couple of antipodal endpoints of some arcs specified, say $P^1$ and $Q^1$.

Denote by $b^1$ the diagrammatic component to which $P^1$ and $Q^1$ belong. Now travel along the boundary of the disk in counterclockwise direction, starting from these two antipodal endpoints. Each time a couple of antipodal endpoints is encountered consider whether it belongs to a new diagrammatic component. If this is the case, denote it by $b^2$ and denote this couple of endpoints by $P^2$ and $Q^2$. Continuing in this way, call the subsequent diagrammatic components from $b^3$ to $b^n$ and, for each one of them, call the couple of antipodal points encountered $P^3$, $Q^3$ to $P^n$, $Q^n$.

Orient each $b^i$ in such a way that $P^i$ becomes the initial point of the arc to which it belongs.
In the definition below, we will use the arc distance between two points $P$ and $X$ where $P$ is an endpoint of an arc. The original definition of arc distance is extended to this case in the following way: consider a point $P'$ in the interior of the arc to which $P$ belongs, which is such that $X$ is not between $P$ and $P'$. Then the arc distance between $P$ and $X$ is by definition the arc distance between $P'$ and $X$.

With these conventions the definition is:

\begin{definition}\label{simple_desc}
The simple diagram $D$ is called {\it descending} with respect to the couple of antipodal points $P^1$ and $Q^1$, provided that for each crossing between $b^i$ and $b^j$, $i<j$, the branch in $b^i$ is over (resp. under) the branch in $b^j$, if the arc distance between $P^i$ and the branch in $b^i$ is even (resp. odd).
\end{definition}

\subsection{General descending diagrams.}\label{gendescdgms}
Consider a diagram $D$. Suppose that $D$ has oriented based 0-homologous diagrammatic components, $a^1, ..., a^m$ and that for each 1-homologous diagrammatic component an {\it ordered} simplifying set is fixed.
Moreover, suppose that a couple of endpoints $P^1$ and $Q^1$ is specified in $D$, that are not in the dashed part of $D$ determined by the crossings in the simplifying sets. With this data we define below the notion of descending diagram.

First we introduce some notations. Let $D'$ be the simple diagram of this data. We will consider $D'$ as a subset of $D$. Note that $P^1$ and $Q^1$ are in $D'$.

Denote by $b^1$ the diagrammatic component to which $P^1$ and $Q^1$ belong. Now travel along the boundary of the disk in counterclockwise direction, starting from these two antipodal endpoints, and consider the successive endpoints encountered. If they belong to a new diagrammatic component {\it and are in the simple diagram $D'$}, denote it by $b^2$ and denote this couple of endpoints by $P^2$ and $Q^2$. Continuing in this way, call the subsequent diagrammatic components from $b^3$ to $b^n$ and, for each one of them, call the couple of antipodal points encountered $P^3$, $Q^3$ to $P^n$, $Q^n$.

Finally for each $b^i$ call its ordered simplifying set $M^i$.

With these conventions the definition is:

\begin{definition}\label{main_def}
$D$ is called {\it descending} with respect to the preceding data, provided that its simple diagram $D'$ is descending with respect to $P^1$ and $Q^1$ and the crossings that are not crossings of $D'$ satisfy the following:

Case 1: Suppose that a crossing $X$ has one branch in $a^i$ and the other either in $a^j$ with $i\le j$ or in $b^l$ for any $l$.
Then, the first pass of $X$ from the basepoint of $a^i$ is an overpass (resp. underpass), if the arc distance between the basepoint and this first pass is even (resp. odd).

Case 2: Suppose that a crossing $Y$ between two 1-homologous diagrammatic components, which is not in a simplifying set, has one branch in a dashed part determined by a crossing $X$ in $M^i$ and for the other branch of $Y$ one of the following holds:

\begin{enumerate}
\item it is in $b^j$ with $j>i$,

\item it is in $b^j$ with $j\le i$, but not in a dashed part of it,

\item it is in $b^i$, in the dashed part determined by a crossing in $M^i$ that is greater or equal to $X$ ($M^i$ is ordered).
\end{enumerate}

Give $b^i$ the orientation determined by $X$.
Then, the first pass of $Y$ from the upper branch of $X$ is an overpass (resp. underpass), if the arc distance between the upper branch of $X$ and this first pass of $Y$ is even (resp. odd).
\end{definition}

\subsection{Unknottedness of simple descending diagrams}
Suppose that in a simple diagram there is a diagrammatic component $a$ which has a unique arc and another diagrammatic component $b$. In a simple diagram all diagrammatic components are 1-homologous, so $a$ divides the disk in two halves. An arc of $b$ can have its endpoints in different halves (in which case we say that it {\it crosses} $a$) or in the same half (in which case we say that it {\it does not cross} $a$). In order to prove that simple descending diagrams are diagrams of standard unlinks, we need the following:

\begin{lemma}\label{cross}
Let $D$ be a simple diagram. Let $a$ and $b$ be two of its diagrammatic components. Suppose that $a$ has a single arc and that $b$ has $n$ arcs, $n\ge 3$.
Then at least two of the arcs of $b$ do not cross $a$.
\begin{proof}
First we prove that there is at least one such arc.

Suppose that all arcs of $b$ cross $a$. Consider the net of $D$, say $N$. Let $P$ be the image in the net of the endpoints of $a$. If we remove from $N$ the image of $a$ and a small disk centered at $P$ we get the Mobius strip $M$ together with some closed curves on it. By assumption the image of $b$ is a simple closed curve in $M$ representing $\pm n$ in the fundamental group of $M$ (the infinite cyclic group). This is impossible if $n\ge 3$ by \cite{DRJC} (p.145).
Thus at least one of the arcs of $b$ does not cross $a$.

Furthermore, the number of crossings between $a$ and $b$, two 1-homologous diagrammatic components, is odd.
We can count the crossings coming from $b$-arcs that do not cross $a$ (even number) and the other ones (odd number). If there was just one $b$-arc that does not cross $a$, the number of the remaining arcs being even ($b$ is 1-homologous) the total number of crossing between $b$ and $a$ would be even.
This is impossible, so there has to be another $b$-arc that does not cross $a$.
\end{proof}
\end{lemma}

An arc of a diagrammatic component $b$ will be called {\it most nested} if it is possible to travel on the boundary circle of the diagram from one endpoint of the arc to the other without encountering endpoints that belong to other arcs of $b$.

If $P$, $Q$ is a couple of antipodal points on the boundary circle of the diagram and $b^1$ and $b^2$ are two arcs of $b$, we say that $b^1$ is {\it nested} in $b^2$ {\it with respect to} $P$, $Q$ if it is possible to travel on the boundary circle from $P$ to $Q$ and meet successively: an endpoint of $b^2$, the endpoints of $b^1$ and finally the other endpoint of $b^2$.

\begin{theorem}\label{main_lemma}
Suppose that $D$ is a simple diagram with a couple of antipodal endpoints $P^1$ and $Q^1$ specified. If $D$ is descending with respect to $P^1$ and $Q^1$, then $D$ is a diagram of standard unlink.

\begin{proof}
Suppose that $D$ is descending with respect to $P^1$ and $Q^1$.

We will consider successively the diagrammatic components $b^1, ..., b^n$ and the couples $P^1$, $Q^1$ to $P^n$, $Q^n$ defined in the same way as in section \ref{simple_desc_dgms}.

First suppose that $b^1$ has a single arc which has $P^1$ and $Q^1$ as endpoints and lies above everything else. Suppose that $b^2$ has $k$ arcs, $k\ge 3$.
Denote them by $b^2_1$, ..., $b^2_k$ in such a way that $P^2$ is an endpoint of $b^2_1$, $Q^2$ is an endpoint of $b^2_k$ and the order comes from some orientation of $b^2$.

We will reduce the number of arcs of $b^2$. We choose the first arc to eliminate in such a way that it is most nested, does not cross $b^1$ and does not have $P^2$ or $Q^2$ as endpoint.

To see that it is always possible to find such arc we consider different cases, depending on whether $b^2_1$ or $b^2_k$ crosses $b^1$ or not. These cases are pictured in Figure \ref{3cases}.

Cases 1 and 2: If at least one of $b^2_1$ and $b^2_k$ crosses $b^1$, we know that there is another arc, say $b^2_i$, that does not cross $b^1$ by Lemma \ref{cross}. If it is most nested we will eliminate it. Otherwise it is easy to see that there is a most nested arc, say $b^2_j$, nested in $b^2_i$ with respect to the endpoints of $b^1$. Then $b^2_j$ does not cross $b^1$, it does not have $P^2$ or $Q^2$ as endpoint and we will eliminate it.

Case 3: From Figure \ref{3cases} it is clear that there is an arc nested in $b^2_k$ with respect to the endpoints of $b^1$. It does not cross $b^1$. As in cases 1 and 2, we will eliminate it or eliminate another most nested arc, nested in it with respect to the endpoints of $b^1$.

We consider also the situation where the roles of $b^2_1$ and $b^2_k$ are reversed in Case 3. We get a suitable arc to eliminate in that case, too.

\begin{figure}[h]
\scalebox{1}{\includegraphics{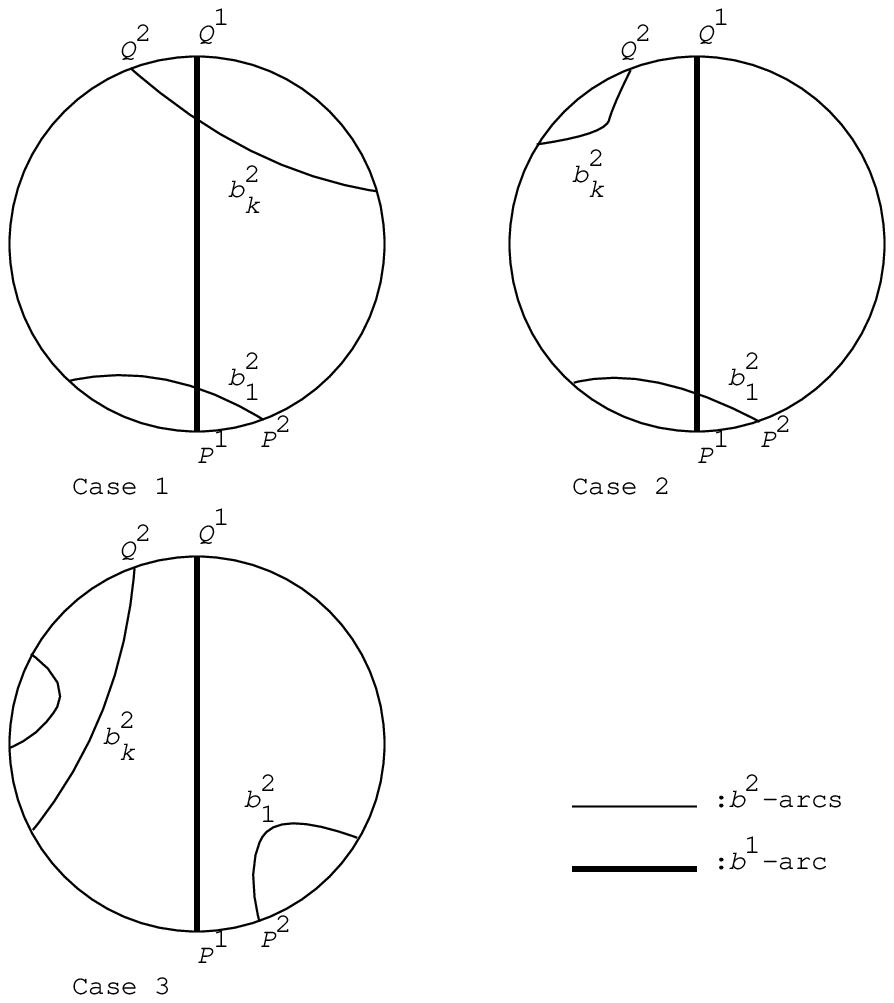}}
\caption{}
\label{3cases}
\end{figure}

To eliminate an arc, consider its position with respect to the arcs of $b^3$, ..., $b^k$. If it lies above these arcs, one of its endpoints is moved towards the other above all endpoints of other arcs. This is done by a sequence of $\Omega_5$ moves. In the next step all crossings between this arc and other arcs are killed with some $\Omega_1-\Omega_3$ moves. Finally it is eliminated with $\Omega_4$.

If an arc lies below the arcs of $b^3$, ..., $b^k$, the same moves are used except that the endpoint of the arc is moved below other endpoints.

Note that it is crucial that the chosen arc does not cross $b^1$. Otherwise, after some application of $\Omega_5$ move between the arc and $b^1$, there would be a crossing between $b^1$ and $b^2$ in which the upper branch would belong to $b^2$. Thus $b^1$ would no more lie above everything else.

After the elimination, a new arc lies above or below $b^3, ..., b^n$, and the number of arcs of $b^2$ is decreased by two.

In this way, the number of arcs of $b^2$ can be reduced to one. At the end $b^1$ and $b^2$, each with a single arc, lie above all other arcs. If we travel from $P^1$, $Q^1$ in counterclockwise direction, the first couple of endpoints encountered is $P^2$ and $Q^2$. This is shown in Figure \ref{2unlink}.

\begin{figure}[h]
\scalebox{1}{\includegraphics{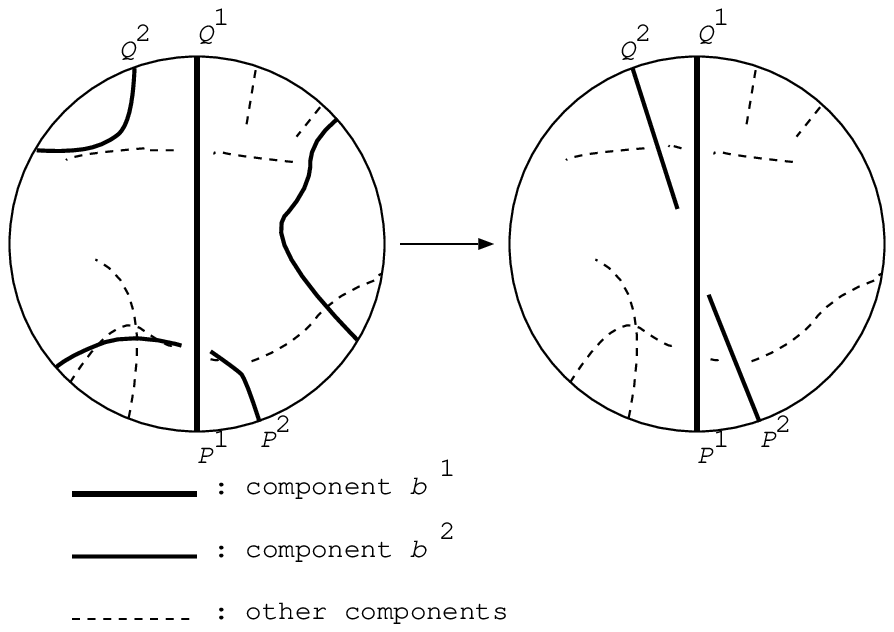}}
\caption{}
\label{2unlink}
\end{figure}

We continue similarly with $b^3$. Its arcs can cross both $b^1$ and $b^2$ or none of them. So the reduction of arcs that was done for $b^2$ works for $b^3$ as well.

In this way we can continue for all arcs and at the end we get standard unlink.

It remains to prove that at the beginning of the isotopy, the number of arcs of $b^1$ can be reduced to one arc which has $P^1$ and $Q^1$ as endpoints and lies above everything else. To see this we add a new diagrammatic component $b^0$ which has a single arc lying above everything else, with $P^0$ and $Q^0$ as endpoints, which are such that if we travel from $P^0$, $Q^0$ in counterclockwise direction, the first couple of endpoints encountered is $P^1$ and $Q^1$.
Now the number of arcs of $b^1$ is reduced to one in the same way as it was done for $b^2$. Finally we remove $b^0$.
\end{proof}
\end{theorem}

\subsection{Unknottedness of general descending diagrams}

\begin{theorem}\label{main_theorem}

Let $D$ be a diagram of a link in $\mathbb RP^3$. Suppose that its 0-homologous diagrammatic components are ordered and each of them is oriented and based. Moreover, suppose that for each 1-homologous diagrammatic component of $D$, an ordered simplifying set is chosen. Finally let $P$ and $Q$ be a couple of antipodal endpoints in the simple diagram of $D$. If $D$ is descending with respect to this data, then $D$ is a diagram of standard unlink.

\begin{proof}
Suppose that $D$ is descending.

Denote the 0-homologous diagrammatic components of $D$ by $a^1, ..., a^m$, according to the order. Denote $P$ and $Q$ by $P^1$ and $Q^1$ respectively. In the proof, we use the same notation as in section \ref{gendescdgms} for the simple diagram of $D$, 1-homologous diagrammatic components, endpoints and simplifying sets.

First, the 0-homologous diagrammatic components are unlinked as in the case of knots, starting with $a^1$ and ending with $a^m$. At the end of this step, we get a diagram in which there are crossings only between 1-homologous diagrammatic components.

Next, consider the first crossing in the ordered simplifying set $M^1$, say $X$. Consider the dashed part determined by $X$. Denote by $b^X_1$, ..., $b^X_k$ the arcs and part of arcs encountered successively when traveling on the net according to the orientation determined by $X$, from the upper branch to the lower branch of $X$. Note that $k$ is odd because of the definition of orientation determined by a self-crossing. As $D$ is descending we have:

$$b^X_2\le b^X_4\; ...\le b^X_{k-1}\le {everything\; else}\le b^X_k\; ...\le b^X_3\le b^X_1$$

Now, as in the proof of Theorem \ref{theorem_knots}, we can eliminate $b^X_2$, $b^X_4$, etc. Finally, we get a single descending part of arc from the upper branch to the lower branch of $X$ that is eliminated at the end. In this way, the dashed part determined by $X$ is erased.

Similarly, we erase all dashed parts determined by crossings in $M^1, ..., M^n$.

Thus, we get the simple diagram of $D$ together with some 0-homologous diagrammatic components. But the simple diagram of $D$ is a diagram of standard unlink because of Theorem \ref{main_lemma}. Thus $D$ is a diagram of standard unlink.
\end{proof}
\end{theorem}

As in the case of knots the last result implies:

\begin{corolary}\label{nonoriented}
Let $D$ be a diagram of a link. By making some crossing changes on it, we can obtain a diagram of standard unlink.
\end{corolary}

\section{The oriented case}\label{counter}
\begin{figure}[h]
\scalebox{1}{\includegraphics{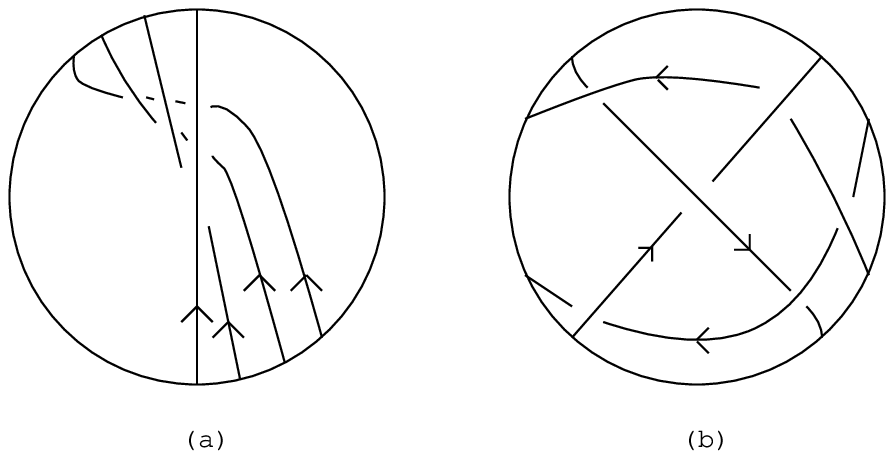}}
\caption{}
\label{4unlink}
\end{figure}

The next natural step is to try to unlink an oriented link. Now, for standard oriented unlink we can take a standard unlink and orient each 1-homologous component in such a way that at each crossing the local writhe is +1.
For instance, standard oriented unlink for four 1-homologous components is pictured in Figure \ref{4unlink} (a).

The notion of descending diagram cannot be extended to the oriented case in a satisfactory way. In other words, there is no way to obtain the equivalent of Corollary \ref{nonoriented} for oriented links: there is a diagram such that no matter what crossing changes are done on it, it will never become a standard oriented unlink.

A counterexample is presented in Figure \ref{4unlink} (b). This link with four 1-homologous components has a local writhe +1 at each crossing.
Each couple of components has a unique common crossing so, if we want to obtain the standard oriented unlink, we cannot change any crossing because we would get at this crossing a local writhe of -1. But the link in Figure \ref{4unlink} (b) is not the standard oriented unlink.
It is the mirror image of the link $6_3^4$, that appears in \cite{JD2}, whereas the standard oriented unlink is the mirror image of the link $6_2^4$, and these two links are not isotopic.

\end{document}